\documentclass[11pt]{article}
\usepackage[margin=2.7cm]{geometry}
\usepackage[T1]{fontenc}
\usepackage{lmodern}
\usepackage{booktabs}
\usepackage{microtype}
\usepackage[hidelinks]{hyperref}
\usepackage{orcidlink}

\title{Improved upper bounds for five IHTC-2024 benchmark instances\\
via exact reoptimization of the nurse-to-room layer}
\author{Alfonso Gippini Requeijo\,\orcidlink{0009-0004-2109-0583}\\
\small Team Banzai S.L., Madrid, Spain\\
\small \texttt{alfonso@team-banzai.com}}
\date{2 September 2026}

\begin{document}
\maketitle

\begin{abstract}
The Integrated Healthcare Timetabling Competition 2024 (IHTC-2024) defined a benchmark that integrates surgical case planning, patient admission scheduling and
nurse-to-room assignment. Its instance set remains live: the organizers keep a
public table of best-known upper bounds that is open to new submissions. We report five improved upper bounds on that table, for instances i02 (1263, previously 1264), i05
(12744, previously 12760), m01 (3380, previously 3384), m03 (6692, previously 6697) and
m04 (3315, previously 3318). All five solutions were checked with the official IHTP
validator and carry no hard-constraint violations. The solutions are not built from
scratch: each is \emph{derived} from the corresponding solution published on the
benchmark page, holding its admission days, room assignments and operating-theater
assignments fixed, and reoptimizing only the nurse-to-room-and-shift layer to proven
optimality. The improvements are small in magnitude, between 0.075\% and 0.125\%, but they
carry a structural observation: on
the instances where the nurse layer can be closed exactly, it closes within seconds on all
but one of them, and several of the published state-of-the-art solutions fall short of
optimality in that layer.
\end{abstract}

\section{Introduction}
The IHTC-2024 competition \cite{ceschia2025ihtc,ihtcweb} posed an integrated hospital
timetabling problem that combines three sub-problems usually studied separately:
surgical case planning, patient admission scheduling, and nurse-to-room assignment. The
competition drew 32 submissions, five of which were chosen as finalists and
evaluated on a hidden dataset.

The competition has ended, but the benchmark has not. The organizers keep a
page of best-known upper bounds for the 30 public instances (i01 to i30) and the 30
hidden instances (m01 to m30) \cite{ihtpbounds}, and explicitly accept new upper bounds
sent as JSON solution files. As of 22 August 2026, that table had already advanced
past the competition results on several instances.

This note records five further improvements to that table. The contribution is
deliberately narrow. We do not propose a new solution method for the integrated problem. We show that a specific, well-understood component of the problem,
the nurse-to-room assignment with everything else held fixed, is still not solved to
optimality in several of the currently published best solutions, and that closing it
exactly is cheap.

\paragraph{Relation to prior work.}
Exact reoptimization of the nurse-to-room layer with the patient-side decisions fixed is not new. Team Twente used precisely this idea as Phase~3 of their competition entry, and
analyzed its contribution separately \cite[Sections~4.4 and~5.5]{guericke2026hybrid}.
The same two-step structure appears in the public repository of another competition
participant \cite{team42repo}. Multi-neighborhood approaches to the related integrated
patient-to-room and nurse-to-patient problem have also been studied by the organizers'
group \cite{zanazzo2025multi}. The novelty here is the target: we apply the exact nurse step \emph{a
posteriori} to the best solutions currently published on the benchmark page, produced by
several different teams with several different methods, and we report what that exposes
about them.

\section{The benchmark and its living upper-bound table}
\label{sec:bench}
\paragraph{Instance set.} The public dataset holds 30 instances (i01 to i30) spanning
horizons of two to four weeks, with up to roughly 500 patients and between 10 and 65
nurses. The hidden dataset used in the second stage holds 30 further instances (m01
to m30). Both sets, together with the winners' solutions, are published by the organizers \cite{ihtcweb}.

\paragraph{Objective and feasibility.} A solution is scored by a weighted sum of
eight soft-cost components, under nine hard constraints. The arbiter is the
official validator shipped with the competition material \cite{validator}, which
reports both the total cost and the number of hard-constraint violations. Every value
reported in this note is the figure printed by that validator, run on the solution file we provide.

\paragraph{Baseline.} Our baseline is the best-known upper bound listed on the living
IHTP page \cite{ihtpbounds}, consulted on 19 August 2026 and re-checked on 22 August
2026, at which date all five baseline values stood unchanged. For i02, i05, m01 and m03
the listed bound matches the value attained during the competition. For m04 the
listed bound (3318) is itself an improvement over the value of the published competition solution (3345), so
our result improves a bound that had already moved after the competition.

\section{Method}
We solve a single restricted subproblem. Given an instance and a feasible published
solution for it, we fix every patient-side decision exactly as it appears in that
solution: the admission day of each scheduled patient, the room assigned for the whole length of stay, and the assigned operating theater. The set of unscheduled optional
patients is likewise untouched. All soft-cost components that depend only on those
decisions are therefore constant, and the residual problem is the assignment of nurses
to room-and-shift pairs. We state that residual problem as a constraint program and
solve it with CP-SAT from Google OR-Tools \cite{ortools}, with a budget of 180 seconds
per instance and four threads. A result is reported only when the solver returns a proven-optimal status for this residual problem. The resulting full solution is then re-evaluated end to end by the official validator, and is accepted only if it is
feasible and its total cost is strictly below the published bound. The residual
model is solved from scratch, without warm-starting from the published nurse assignment
and without a cost-improvement cut, so a non-proof tells us nothing about the
published solution.

\section{Results}
We ran the procedure on all 60 published solutions of the table. The residual
nurse problem was closed to proven optimality on 9 of the 60 instances. On those 9, the
published solution was already optimal in the nurse layer in 4 cases (i01, i03, i06,
m02), and improvable in the remaining 5. On the other 51 instances the solver did not
prove optimality within the budget. All nine instances on which the layer
closed sit at the small end of the set (28 to 105 patients, at most 10 rooms); none of
the 46 instances with more than 105 patients closed, so the observation below
holds solely in that regime.

Table~\ref{tab:results} collects the five improvements. Lower bounds are taken from
\cite[Table~5]{guericke2026hybrid}, where they were obtained from one-hour runs of an
independent method; they are included to frame the remaining gap.

\begin{table}[ht]
\centering
\caption{Improved upper bounds for five IHTC-2024 instances. ``Previous UB'' is the
best-known upper bound listed on the living benchmark page \cite{ihtpbounds} as of 19
August 2026 and unchanged on 22 August 2026. ``New UB'' is the total cost reported by the
official validator \cite{validator} on the solution file we submit, which in all five
cases shows zero hard-constraint violations. ``LB'' is the lower bound published in
\cite[Table~5]{guericke2026hybrid}. ``Time'' is the wall-clock time in which the
residual nurse problem was closed to optimality. Gap to LB is computed as (UB $-$ LB)/LB.}
\label{tab:results}
\small
\begin{tabular}{lrrrrrrrr}
\toprule
Instance & Prev.\ UB & New UB & Impr. & Impr.\ [\%] & LB & Gap to LB [\%] & Time [s] & Viol. \\
\midrule
i02 & 1264  & \textbf{1263}  & 1  & 0.079 & 1049  & 20.40 & 0.5   & 0 \\
i05 & 12760 & \textbf{12744} & 16 & 0.125 & 12642 & 0.81  & 6.5   & 0 \\
m01 & 3384  & \textbf{3380}  & 4  & 0.118 & 2968  & 13.88 & 1.2   & 0 \\
m03 & 6697  & \textbf{6692}  & 5  & 0.075 & 6419  & 4.25  & 0.4   & 0 \\
m04 & 3318  & \textbf{3315}  & 3  & 0.090 & 2544  & 30.31 & 154.0 & 0 \\
\bottomrule
\end{tabular}
\end{table}

\paragraph{Experimental conditions.} Runs were carried out on a single desktop machine (AMD
Ryzen~7 7800X3D, 8 physical cores, 64\,GB RAM) under Manjaro Linux, with OR-Tools
9.15.6755 on Python 3.14, pinned to four physical cores. The machine was \emph{not}
dedicated: a separate workload of the author ran alongside the
sweep, and CPU contention was present throughout. A proof of optimality does not depend on
how much CPU was free (m04, at 154\,s of a 180\,s budget, is the single result that a busier
machine might have lost). A feasible solution of
lower cost is likewise unaffected, since the validator does not ask what hardware produced it. Contention can only have affected the negative side, the 51 instances without a proof.

This budget is not comparable with the competition protocol, which allowed 10 minutes and
four threads per instance for solving the full integrated problem from scratch
\cite{ihtcweb}. Our procedure starts from an already-published solution and optimizes
one layer of it. It is a post-processing step, not a solver, and the numbers in
Table~\ref{tab:results} should not be read as a competitive result against the
finalists.

\paragraph{Status of the benchmark page.} The five solution files went to the
organizers on 24 August 2026. On 25 August 2026 the benchmark site swapped the five files
linked from the living table for ours: each file the table now links is byte-identical to the
corresponding file in our deposit (Section~\ref{sec:data}), whose validator costs are those of
Table~\ref{tab:results}. The same update moved four cells of the cost column, from 1264 to
1263, from 12760 to 12744, from 3384 to 3380 and from 6697 to 6692. The fifth cell was left
alone: for m04 the linked file is ours and evaluates to 3315, while the column showed 3318 when
we last consulted the page, on 27 August 2026. The mismatch is in that cell, not in the file it
links to, and we have told the organizers; it may well have been fixed by the time this note is
read. The date and the extent of
the update can be checked in the public repository of the site \cite{ihtcrepo}. The
``Previous UB'' column of Table~\ref{tab:results} accordingly refers to the state of the page
before that update.

\section{Discussion}
The magnitudes in Table~\ref{tab:results} are between 1 and 16 cost units.

First, the residual nurse problem is easy where it is tractable at all. Three of the five
improvements were proven optimal in under 1.5 seconds, and one in 6.5 seconds. The
outlier, m04, took 154 seconds. On the instances where the layer closes, it closes almost immediately.

Second, and in spite of this, five of the nine published solutions on which the layer could
be closed were suboptimal in it. These are solutions produced by strong methods,
including competition finalists, and they are the current best-known values on a
maintained benchmark. That the cheapest component of the problem is still leaving cost
on the table in more than half of the cases where it can be checked points to a
nurse layer that is systematically under-exploited in the final polishing stage, even by
approaches that treat it explicitly.

Third, the remaining gaps to the published lower bounds (Table~\ref{tab:results},
ranging from 0.81\% on i05 to 30.31\% on m04) confirm that the dominant difficulty of
this benchmark lies in the patient-side decisions. Closing the
nurse layer exactly does not approach the lower bound; it only removes a small, reliably
recoverable residue.

\paragraph{Limitations.} The improvements are derivative: they inherit the patient-side
decisions of solutions produced by other teams, and the credit for most of the
objective value belongs to whoever produced those base solutions. The benchmark page
does not attribute its listed bounds to specific teams, so we cannot name them
here, and we would welcome the organizers' guidance on the appropriate attribution.
Finally, our sweep used a 180-second budget; a longer one would likely close the
nurse layer on more than 9 of the 60 instances.

\section{Data availability}
\label{sec:data}
The five solution files are archived on Zenodo under DOI
\href{https://doi.org/10.5281/zenodo.22061899}{10.5281/zenodo.22061899} (v1; all versions: 10.5281/zenodo.22061898) under a CC-BY 4.0 license,
in the JSON format set by the competition, together with a pipeline certificate and
the SHA-256 digest of every file, so that any reader can reproduce the reported costs by
passing them through the official validator \cite{validator} against the corresponding
official instance files. The source code that produced the solutions is not released.

\section*{Acknowledgments}
We thank the organizers of the IHTC-2024, Sara Ceschia, Roberto Maria Rosati, Andrea
Schaerf, Pieter Smet, Greet Vanden Berghe and Eugenia Zanazzo, for designing and
publishing the benchmark, for distributing the official validator, and above all for
maintaining the living table of best-known bounds after the competition ended, without
which none of the results in this note would have been measurable. We also thank the
competition participants whose published solutions form the basis of the five solutions
reported here.

Two debts are more specific. Daniela Guericke, Rolf van der Hulst, Asal Karimpour, Ieke Schrader
and Matthias Walter published the lower bounds that fill the last columns of
Table~\ref{tab:results} \cite{guericke2026hybrid}; the exact nurse step applied here is the same idea they had used as a phase of their own
solver. The work of Eugenia Zanazzo, Sara Ceschia and Andrea Schaerf on the related
patient-to-room and nurse-to-patient problem \cite{zanazzo2025multi}, and the entry that Team~42
released openly \cite{team42repo}, both helped us understand where this step fits.

We thank the free-software community as well. All the computation behind this note ran on CP-SAT
from Google OR-Tools, maintained by Laurent Perron and Vincent Furnon \cite{ortools}.

\section*{Use of AI tools}
This note was written by the author with substantial assistance from a large language model. The
author's input was dictated by voice and transcribed automatically; the model produced and
translated much of the English prose, which the author then directed, edited and checked against
the underlying data. The author takes full responsibility for the contents.

The results are not the output of a language model. Every number reported here comes from the
solver runs and the certificates described in the Method section, and no figure was written into this note
without being checked against the artifact that produced it.

\section*{Competing interests}
The author is the director of Team Banzai S.L., a company that builds optimization
software. No external funding was received for this work.

\end{document}